\input amstex
\documentstyle{gen-j}
\input epsf.tex
\chardef\ss="19
\def\3{\ss}
\def\<{<\!\!<}

\def\Ad{\mathop{\text{\rm Ad}}\nolimits}

\def\conv{\mathop{\text{\rm conv}}\nolimits}

\def\eff{\mathop{\text{\rm eff}}\nolimits}

\def\fix{\mathop{\text{\rm fix}}\nolimits}
\def\Fix{\mathop{\text{\rm Fix}}\nolimits}

\def\id{\mathop{\text{\rm{id}}}\nolimits} 
\def\im{\mathop{\text{\rm{im}}}\nolimits}

\def\Int{\mathop{\text{\rm int}}\nolimits}

\def\pr{\mathop{\text{\bf pr}}\nolimits}

\def\reg{\mathop{\text{\rm reg}}\nolimits}

\def\Sp{\mathop{\text{\rm Sp}}\nolimits}

\def\tr{\mathop{\text{\rm tr}}\nolimits}

\def\0{{\bf 0}}
\def\1{\text{\bf {1}}}

\def\a{{\frak a}}

\def\b{{\frak b}}

\def\g{{\frak g}}

\def\k{{\frak k}}

\def\m{{\frak m}}
\def\n{{\frak n}}

\def\p{{\frak p}}
\def\q{\mathop{\text{\bf q}}}

\def\s{{\frak s}}

\def\so{{\frak {so}}}

\def\t{{\frak t}}
\def\uu{{\frak u}}

\def\v{\mathop{\text{\bf v}}}

\def\x{\mathop{\text{\bf x}}}
\def\y{\mathop{\text{\bf y}}}
\def\z{{\frak z}}

\def\B{\mathop{\text{\bf B}}\nolimits}
\def\C{{\Bbb C}}

\def\R{{\Bbb R}}

\def\:{\colon}  
\def\.{{\cdot}}
\def\|{\Vert}

\def\giantskip{\vskip2\bigskipamount}

\def \la {\langle}
\def\msk{\medskip}
\def \ra {\rangle}
\def \res {\!\mid\!\!}

\def\ssk{\smallskip}

\def\giantbreak{\par \ifdim\lastskip<2\bigskipamount \removelastskip
         \penalty-400 \giantskip\fi}

\def\nin{\noindent}
\def\pagebreak{\vskip 0pt plus 0.0001fil\break}
\def\linebreak{\break}

\def\nin{\noindent}
\def\oline{\overline}

\def\phi{\varphi}


\def\subeq{\subseteq}
\def\supeq{\supseteq}

\def\tilde{\widetilde}

\def\bs{\backslash}

\def\SO{\mathop{\text{\rm SO}}\nolimits}

\topmatter 
\title Lagrangian submanifolds and moment convexity  \endtitle

\rightheadtext{Lagrangians and moment convexity}
\leftheadtext{B. Kr\"otz and M. Otto}

\author Bernhard Kr\"otz  and Michael Otto\endauthor 

\abstract
We consider a Hamiltonian torus action $T\times M \rightarrow
M$ on a compact connected symplectic manifold $M$ and its 
associated  momentum map $\Phi$. For certain Lagrangian submanifolds 
$Q\subseteq M$ we show that $\Phi(Q)$ is convex. 
The submanifolds $Q$ arise as the fixed point set of an involutive
diffeomorphism $\tau:M\rightarrow M$ which satisfies several
compatibility conditions with the torus action, but which is in
general not anti-symplectic. As an application we complete a
symplectic proof of Kostant's nonlinear convexity theorem.
\endabstract

\address Department of Mathematics, University of Oregon, 
Eugene OR 97403-1221 \endaddress 

\address The Ohio State University, Department of Mathematics, 
231 West 18th Avenue, Columbus OH 43210-1174\endaddress 
\email kroetz\@math.uoregon.edu, otto\@math.ohio-state.edu\endemail

\thanks The work of the first author was supported in part by NSF-grant DMS-0097314\endthanks

\endtopmatter 
\document 

\head 1. Introduction\endhead 

In the context of Hamiltonian torus actions $T\times M\to M$ 
on a connected symplectic manifold $M$  one is interested in 
convexity properties of the image of the associated momentum 
map $\Phi\: M\to \t^*$.  This is because of its many 
applications to classical eigenvalue problems and their 
Lie theoretic generalizations. In this paper we will
determine a class of Lagrangian submanifolds $Q\subeq M$
for which $\Phi(Q)$ is convex. Applications to
a symplectic proof of Kostant's  non-linear 
convexity theorem will be given.

\par Before we will describe our results in more detail, it is 
useful to summarize  some known convexity results for the momentum map
$\Phi$. We recall that the critical set of $\Phi$ is the
set $\Fix(M)$ of $T$-fixed points in $M$. Then the convexity 
theorem of Atiyah-Guillemin-Sternberg [2,5] reads as follows.  

\proclaim {Theorem 1.1} If $M$ is compact, then $\Phi(M)$ is convex. 
More precisely, $\Phi(M)$ is the convex polyhedron spanned by the finite 
set $\Phi(\Fix(M))$.\endproclaim 

This theorem has been generalized by Duistermaat [4].  
Assume that $M$ carries an anti-symplectic 
involution $\tau\: M\to M$ 
and write $Q$ for the fixed point set of $\tau$. 
We require  that $Q$ is non-empty. 
Then $Q$ is a Lagrangian submanifold of $M$, and Duistermaat's
convexity theorem [4, Th. 2.5] says: 

\proclaim {Theorem 1.2} If $M$ is compact,
and $\tau\: M\to M$ is an anti-symplectic involution 
which satisfies $\Phi\circ \tau=\Phi$,
then $\Phi (Q)=\Phi(M)$. 
In particular $\Phi(Q)$ is convex. \endproclaim 

Duistermaat's Theorem was further generalized 
to the case where $M$ is non-compact and the momentum map proper 
[9, 14].  
But even if $\Phi$ is not proper, there are 
interesting classes of symplectic 
manifolds for which $\Phi(M)$ is still convex (cf.\ [7,
10,14]). 

\par In this paper we want to give another generalization of
Theorem 1.2 which goes in a different direction. It turns out 
that the assumption that the involution $\tau$ is anti-symplectic
is too strong for certain applications. However, we will show 
that under relaxed conditions on $\tau$ we can still derive the
results of Duistermaat's theorem. In particular,

\proclaim {Theorem 3.1} Let $M$ be a compact connected 
symplectic manifold with Hamiltonian torus action 
$T\times M\to M$ and momentum map $\Phi\: M\to \t^*$. 
In addition, let $\tau\: M\to M$ be an involutive 
diffeomorphism with fixed point set $Q$ such that 
\roster
\item $t\circ \tau=\tau\circ t^{-1}$ for all $t\in T$. 
\item $\Phi\circ \tau=\Phi$. 
\item $Q$ is a Lagrangian submanifold of $M$. 
\endroster 
Then $\Phi(Q)=\Phi(M)$. In particular 
$\Phi(Q)$  is a convex subset of $\t^*$.
Moreover, the same assertions hold if $Q$ is replaced with
any of its connected components.
\endproclaim

As our main application of Theorem 3.1 we will
complete the symplectic proof 
of Kostant's non-linear convexity theorem 
as given in [13].
\par Let us briefly recall Kostant's result.
Let $G=NAK$ be an Iwasawa
decomposition of a semisimple linear Lie 
group $G$. Write $\tilde a\: G\to A$ for the 
associated middle projection. The Lie algebras 
of $G, N, A$ and $K$ shall be denoted by $\g, \n, \a$
and $\k$. Then Kostant's Theorem [11] asserts
$$(\forall X\in \a)\qquad 
\log \tilde a(K\exp(X))=\conv ({\Cal W}.X)\leqno(1.1)$$
where $\conv({\Cal W}.X)\subeq \a$ denotes the 
convex hull of the Weyl group orbit ${\Cal W}.X$. 
\par Lu and Ratiu [13] were able to deduce Kostant's result
from the AGS-convexity theorem for a complex group $G$.
\par If $G$ is not complex the situation is different. For those 
groups for which $\m=\z_\k(\a)$ is abelian one can show (1.1)
using Duistermaat's Theorem [13], since in these cases the 
involution $\tau$ one encounters is indeed anti-symplectic 
[8].
\par If $\m$ is not abelian, $\tau$ still satisfies the
assumptions 1.-3. in Theorem 3.1. Therefore, Theorem 3.1 can
be used to give a symplectic proof of Kostant's theorem for
an arbitrary $G$.

\msk It is our pleasure to thank Robert J. Stanton 
for his very useful advice on structure and presentation 
of the underlying paper. We would also like to thank the referee for his careful work.

\head 2. Local results \endhead 

This section lays the foundation for the proof of the convexity
theorem in Chapter 3. We will give local descriptions for the momentum 
map $\Phi$ and its restriction $\Phi\res_Q$ in Subsection 2.2.
In 2.1, we fix the notation and prove a lemma on the characterization 
of anti-symplectic involutions on a symplectic vector space which is 
needed in 2.2.

\subhead 2.1. Background \endsubhead 

Let $(M, \omega)$ denote a connected 
symplectic manifold with $\dim M=2n$. Then every smooth function 
$f\in C^\infty(M)$ determines a Hamiltonian vector field 
${\Cal X}_f$ on $M$ which is defined by $df=i({\Cal X}_f)\omega$. 
From $\omega$ one obtains the usual Poisson structure om $M$:
$$\{\cdot, \cdot\}\: C^\infty(M)\times C^\infty(M)\to C^\infty (M), 
\ \ (f,g)\mapsto \{ f,g\}=\omega({\Cal X}_f, {\Cal X}_g).$$
\par  Our next datum is a torus $T$ which we require to act 
symplectically on $M$.  Write 
$$ T\times M\to M, \ \ (t, m)\mapsto t.m$$
for this action. 
\par Let $\t$ denote the Lie algebra of $T$. For $X\in \t$ let  $\tilde X$ be  the corresponding 
vector field on $M$, i.e., 
$$\tilde X_m={d\over dt}\Big|_{t=0} \exp(tX).m\qquad (m\in M).$$

\par We will always assume that the action of $T$ on $M$ is Hamiltonian, that 
is there exists a $T$-equivariant Lie algebra homomorphism 

$$\t \to (C^\infty (M), \{\cdot ,\cdot \}), \ \ X\mapsto \Phi_X$$
such that 
$$i(\tilde X)\omega=d\Phi_X \leqno(2.1.1)$$
holds for all $X\in \t$. 
If $\t^*$ denotes the dual of $\t$, then the assignment
$$\Phi\: M\to \t^*; \ \ \la \Phi(m), X\ra =\Phi_X(m) \qquad (m\in M, X\in \t)$$
defines a smooth map, called the {\it momentum map}.

\par Let $\tau\: M\to M$ be an involutive diffeomorphism. We will 
denote by $Q$ its fixed point set, i.e.
$$Q=\{m\in M \: \tau(m)=m\}, $$
and require $Q$ to be non-empty. 
Notice that $Q$ is a closed submanifold
of $M$. If in addition $\tau$ is  anti-symplectic, then 
$Q$ is a Lagrangian submanifold of $M$. 

\par Write $\Fix(M)$ for the set of $T$-fixed points in $M$. 
Equivalently $\Fix(M)$ is the critical set of the momentum 
map $\Phi\: M\to \t^*$. 

\par We will be interested in involutions $\tau$ satisfying
the following conditions:
\msk 
$$\leqalignno {& t\circ \tau=\tau\circ t^{-1}\ \hbox{for all 
$t\in T$.} & (2.1.2)\cr 
& \Phi\circ \tau=\Phi. &(2.1.3)\cr  
& \hbox{$Q$ is a Lagrangian submanifold of $M$.} & (2.1.4)\cr}$$ 
\msk 

Some remarks on (2.1.2)-(2.1.4) are appropriate.

\remark{Remark 2.1.1} (a) Notice that we do not assume that the 
involution $\tau$ is anti-symplectic; however we require 
the fixed point manifold $Q$ to be Lagrangian. 
\par\nin (b) The conditions in Duistermaat's Theorem are stronger 
than (2.1.2)-(2.1.4). In fact, if $\tau$ is anti-symplectic, 
then (2.1.2) and (2.1.3) are equivalent. 
\endremark

Note that condition (2.1.2) is equivalent to

$$d\tau(m)\tilde X_m=-\tilde X_{\tau(m)}\qquad 
(m\in M, \ X\in \t).\leqno(2.1.2a)$$
In particular we have 
$$d\tau(m)\tilde X_m=-\tilde X_m\qquad 
(m\in Q, \ X\in \t).\leqno(2.1.5)$$

\par First we shall investigate conditions
(2.1.2) and (2.1.4) in the linear case, i.e. when 
$(M,\omega)$ is a symplectic vector space with 
linear torus action and linear involution $\tau$. 

\par Let  $(V, \Omega)$ denote 
a finite dimensional symplectic vector space and 
$\tau\: V\to V$ a linear 
involution. As $\tau\circ \tau=\id_V$, the linear 
operator $\tau$ is semisimple with  eigenvalues $+1$ and $-1$. 
Accordingly we have an eigenspace decomposition $V=V_1\oplus V_{-1}$. 

\par Next we endow $(V, \Omega)$ with a linear symplectic
torus action $T\times V\to V$. Then $V$ decomposes 
into {\it fixed} and {\it effective} part 
$$V=V_{\fix}\oplus V_{\eff}\ , \leqno(2.1.6)$$
where 
$$\eqalign{V_{\fix} & =\{ v\in V\: (\forall X\in \t) \ X.v=0\}\cr 
&= \{ v\in V\: (\forall t\in T) \ t.v=v\}\ , \cr }$$
and 
$$V_{\eff}=\t.V.$$
Notice that $\Omega$ is non-degenerate 
when restricted to $V_{\fix}$ or $V_{\eff}$. 
Hence both $V_{\fix}$ and $V_{\eff}$ become
symplectic subspaces of $V$. 

\par The following lemma might be known to many; different versions
of it are frequently encountered in the literature. 
Nevertheless we wish to provide its simple proof. Note that the 
first three statements do not rely on the existence of the torus 
action at all.

\proclaim{Lemma 2.1.2}
Let $(V,\Omega)$ be a symplectic vector space and 
$\tau\: V\to V$ be a linear involution with eigenspace 
decomposition $V=V_1\oplus V_{-1}$. Then the following
statements are equivalent:
\roster 
\item $\tau$ is anti-symplectic. 
\item $V_{1}$ and $V_{-1}$ are Lagrangian. 
\item $V_{1}$ is Lagrangian and there exists a symplectic isomorphism
$\varphi:V\rightarrow V $ with $\varphi(V_{1})=V_{-1}$ and 
$\varphi(V_{-1})=V_{1}$.
\item $V_{1}$ is Lagrangian and there exists a linear symplectic torus 
action $T\times V \to V$ with the following properties: 
\itemitem $t\circ \tau = \tau \circ t^{-1} \quad \forall \  t\in T .$ 
\itemitem $V_{\fix}=\{0\}$.
\endroster 
\endproclaim

\demo{Proof} For  $v\in V$ let  $v_{1}\in V_{1}, v_{-1}\in V_{-1}$be such that
$v=v_{1}+v_{-1}$. 
\ssk\nin 1. $\Longrightarrow$ 2.:  Assume that $\tau$ is anti-symplectic. Then
$$\Omega(v_{1},w_{1})=\Omega\left(\tau (v_{1}),\tau (w_{1})\right)=-\Omega
(v_{1},w_{1}),$$
and similarly,
$$\Omega(v_{-1},w_{-1})=\Omega\left(-\tau (v_{-1}),-\tau (w_{-1})\right)=-\Omega
(v_{-1},w_{-1}). $$
This implies that $V_{1}$ and $V_{-1}$ are isotropic. As $V=V_1\oplus V_{-1}$, 
it follows that both $V_1$ and $V_{-1}$ are maximally isotropic, i.e.
Lagrangian subspaces of $V$. 

\ssk\nin 2. $\Longrightarrow$ 1.:
Assume that  $V_{1}$ and $V_{-1}$ are Lagrangian. We compute

$$\eqalign{\Omega(v_{1}+v_{-1},w_{1}+w_{-1}) & =  \Omega(v_{1},w_{1})+
\Omega(v_{-1},w_{1})+\Omega(v_{1},w_{-1})+\Omega(v_{-1},w_{-1}) \cr 
& = \Omega(v_{-1},w_{1})+\Omega(v_{1},w_{-1}), \cr}$$
and 
$$\eqalign{\Omega(\tau(v_{1}+v_{-1}),\tau(w_{1}&+w_{-1})) = 
\Omega(v_{1}-v_{-1},w_{1}-w_{-1})\cr 
 & =  \Omega(v_{1},w_{1})-
\Omega(v_{-1},w_{1})-\Omega(v_{1},w_{-1})+\Omega(v_{-1},w_{-1}) \cr 
& =  -\Omega(v_{1}+v_{-1},w_{1}+w_{-1}).\cr}$$
Hence, $\tau$ is anti-symplectic. 

\ssk\nin 2. $\Longrightarrow$ 3.:
As Lagrangian subspaces both $V_{1}$ and $V_{-1}$ have dimension
$n=\frac{1}{2} \dim V $. Let $\{ e_1,\dots e_n \}$ be a basis
of $V_1$. As $\Omega$ is non-degenerate, assumption (2) implies that there is 
a basis $\{ f_1, \dots 
f_n \}$ of $V_{-1}$ such that $\Omega(e_i, f_j)=\delta_{ij}$ for all
$1\leq i, j\leq n$. 
\par Define $\varphi$ by 
$$\varphi(e_i)=f_i, \quad \varphi(f_i)=-e_i \quad \forall \  i.$$
Then for all $1\leq i, j\leq n$
$$\Omega(e_i, f_j)=\delta_{ij}=\Omega(e_j, f_i)=-\Omega(f_i, e_j)
=\Omega(\phi(e_i), \phi(f_j)), $$
completing the proof of 1. $\Longrightarrow$ 3.

\ssk\nin 3. $\Longrightarrow$ 2.:
It suffices to show that $V_{-1}$ is isotropic. But this follows
from the fact that $\varphi$ is surjective and that for $v,w\in V_{1}$,
$$ \Omega(\varphi(v),\varphi(w))=\Omega(v,w)=0. $$

\ssk\nin 3. $\Longrightarrow$ 4.:
For $\t=\R \varphi$, the corresponding torus action by $T=\exp \t $ clearly 
has the desired properties.

\ssk\nin 4. $\Longrightarrow$ 2.:
We have to show that $V_{-1}$ is Lagrangian. 
Let $T\times V \to V$ be a torus action
which satisfies (5) and (6). We notice that 
(5) is equivalent to its infinitesimal 
version 
$$\tau\circ X= -X\circ \tau \qquad (X\in \t)\ . 
 \leqno(2.1.7)$$
Let $X\in \t$. Then it is immediate from (2.1.7) that 
$$X(V_1)\subseteq V_{-1}\quad\hbox{and}\quad 
X(V_{-1})\subseteq V_1\ . \leqno(2.1.8)$$
It follows from (6) that 
there is an element $Y\in \t$ such that
$Y\: V\to V$ is invertible. Hence (2.1.8)
implies that 
$$Y(V_1)=  V_{-1}\quad\hbox{and}\quad 
Y(V_{-1}) =V_1\ . \leqno(2.1.9)$$

\par Let $v_{-1},w_{-1}\in V_{-1}$. By (2.1.9) we find 
$v_1,w_1\in V_1$ such that $Y.v_1=v_{-1}$ and 
$Y.w_1=w_{-1}$. 
As $V_1$ is Lagrangian, $Y^2(V_1)=V_1$ (by (2.1.9))
and $\Omega$ is $T$-invariant, it now 
follows that 
$$\Omega(v_{-1},w_{-1})=\Omega(Y.v_1, Y.w_1)
=-\Omega(Y^2.v_1, w_1)=0\ .$$
Hence $V_{-1}$ is Lagrangian. 
\qed \enddemo

\subhead 2.2. Local normal forms\endsubhead 

Throughout this subsection we will assume that 
$(M,\omega)$ is a connected symplectic manifold 
endowed with a Hamiltonian torus action with momentum map $\Phi$.
Also, we have an involution $\tau$ which satisfies 
(2.1.2) - (2.1.4) and whose fixed point set we denote by $Q$.  
Our objective  is to provide a local normal form for $\Phi\res_Q$ 
near a point $m\in Q$.
To that end we first recall a method of finding suitable local descriptions
for $\omega$ and $\Phi$ in the neighborhood of a generic point
$m\in M$. We then consider points $m\in Q$ and obtain a refined form
of $\Phi\res_Q$ which is adapted to the involution $\tau$. 
We start with a simple observation 
(cf.\ [4, Lemma 2.1]):

\proclaim {Lemma 2.2.1} Let $m\in Q$ and $X\in \t$. Then 
$d(\Phi_X\res_Q)(m)=0$ implies $d\Phi_X(m)=0$. In particular, 
$m$ is fixed under the action of the one parameter 
subgroup \linebreak $\exp(\R X)$. \endproclaim 

\demo{Proof} Write $E={\Cal T}_m M$ for the tangent space 
at $m$. Let $E=E_1\oplus E_{-1}$ be the decomposition 
of $E$ into $\pm 1$-eigenspaces of the involution 
$d\tau(m)$. Notice that 
${\Cal T}_mQ=E_1$. In order to show that $d\Phi_X(m)=0$ 
it hence suffices to prove $d\Phi_X(m)(v)=0$ for all $v\in E_{-1}$. 
\par Let $v\in E_{-1}$. Then it follows from (2.1.3) that 
$$d\Phi_X(m)(v)= -d\Phi_X(m)(d\tau(m)v) =-d(\Phi_X\circ \tau)(m)(v)
=-d\Phi_X(m)(v), $$
and so $d\Phi_X(m)(v)=0$. The last assertion in the Lemma
follows from (2.1.1). 
This concludes the proof of the lemma.
\qed \enddemo

For $m\in M$ we write $T_m$ for the stabilizer 
of $T$ in $m$, i.e.
$$T_m=\{ t\in T\: t.m=m\}.$$
The Lie algebra $\t_m$ of $T_m$ is then given by 
$$\t_m=\{ X\in \t\: \tilde X_m=0\}.$$
If in addition $m\in Q$, then it follows from 
Lemma 2.2.1 that we can equally characterize $\t_m$
by 
$$\t_m=\{ X\in \t\: d(\Phi_X\res_Q)(m)=0\}.\leqno
(2.2.1) $$

\msk Fix now $m\in M$. Next we provide local 
normal forms for $\omega$ and $\Phi$ near $m$.
We will recall the procedure of momentum 
reconstruction (cf.\ [6, Ch.\ 41]): The momentum map near $m$ is 
uniquely characterized by $\Phi(m)$, the stabilizer $T_m$ and the 
linear representation of $T_m$ on the tangent space 
${\Cal T}_m M$. This is even true for a general 
compact Lie group $T$ and in case of a torus 
was further  exploited in [9, Sect. 2]. 

\par Notice that ${\Cal T}_m(T.m)$ is an isotropic subspace of ${\Cal T}_mM$. Thus $\omega_m$ 
induces on the quotient 
$$V= {\Cal T}_m (T.m)^\perp/ {\Cal T}_m (T.m)$$
a symplectic form $\Omega$.  We write $\Sp(V, \Omega)$
for the corresponding symplectic group. 
Notice that the isotropy subgroup $T_m$ acts on ${\Cal T}_mM$ symplectically. 
Clearly this action leaves ${\Cal T}_m (T.m)$ and hence ${\Cal T}_m (T.m)^\perp$
invariant, thus giving rise to a representation on $V$, say 
$$\pi\: T_m \to \Sp(V,\Omega).$$

\par Write $(T_m)_0$ for the connected component of $T_m$ 
containing $\1$. Notice that $(T_m)_0<T$ is a subtorus and 
so we can find a torus complement $S_m$ to $(T_m)_0$ in $T$, i.e.
$$T=S_m\times (T_m)_0.$$
We denote the Lie algebra of $S_m$ by $\s_m$. Then 
$\t=\t_m\oplus \s_m$ and we have a canonical identification 
$\t^*=\s_m^*\times\t_m^*$. Denote by ${\Cal T}^*S_m$ the cotangent bundle of 
$S_m$ with its canonical symplectic structure. In the sequel 
we use the identification  ${\Cal T}^*S_m=\s_m^*\times S_m$. 
Hence ${\Cal T}^*S_m\times V$ carries a natural 
symplectic structure. Further $T=S_m\times (T_m)_0$ acts 
symplectically on ${\Cal T}^*S_m\times V$ via 

$$(s,t).(\beta,s' v)\:= (\beta, ss', \pi(t)v)
\leqno (2.2.2)$$
for $s,s'\in S_m$, $t\in (T_m)_0$, $\beta\in\s_m^*$ 
and $v\in V$. 
It follows from [9, Lemma 2.1] that 
there is a symplectic diffeomorphism 
$$\rho\: {\Cal T}^*S_m\times V\supeq {\Cal U}\to U\subeq M$$
from an open neighborhood ${\Cal U}$ of 
$(0,\1,0)\in {\Cal T}^*S_m\times V$
to an open neighborhood $U$ of $m$ such that $\rho$ 
is locally $T$-equivariant and satisfies 
$\rho(0,\1,0)=m$. 

\par In the following we will  identify $U$ with ${\Cal U}\subeq {\Cal T}^*S_m\times V$
via our symplectic, locally $T$-equivariant chart $\rho\: {\Cal U}\to U$. 
Write $V=V_{\fix}\oplus V_{\eff}$ for 
the decomposition of $V$ in effective and fixed part 
for the linear action of $T_m$ on $V$ (cf.\ (2.1.6)).
Furthermore we have the $\t_m$-weight space 
decomposition $V=\oplus_{\lambda\in \Lambda} V_\lambda$. 
We decompose elements $v\in V$ 
as $v=\sum_{\lambda\in \Lambda} v_\lambda$
with $v_\lambda\in V_\lambda$. 
Recall that there is a $T_m$-invariant complex structure $J$ on 
$V_{\eff}$  such that $\la v,w\ra =\Omega(v,Jw)$  
defines a positive definite scalar product on 
$V_{\eff}$. 
Then it follows from [9, Lemma 2.2] that the local normal form of
$\Phi$ near a generic point $m\in M$ is given on ${\Cal U}$ 
by 

$$\Phi\: {\Cal U}\to \t^*=\s_m^*\times \t_m^*, 
\ \ (\beta,s, v)\mapsto \Phi(0,\1,0)+ 
\Big(\beta, {1\over 2}\sum_{\lambda\in\Lambda\atop\lambda\neq 0} 
\|v_\lambda\|^2 \lambda \Big).\leqno(2.2.3)$$

\ssk  Assume now that $m\in Q$. 
So far we have not adressed the question of the
nature of $Q$ and $\tau$  within our  new coordinates 
in ${\Cal T}^*S_m\times V$. In case $\tau$ 
is anti-symplectic on $M$, there is 
a beautiful answer, namely $\tau(\beta, s, v)=(\beta, s^{-1}, \tau_V(v))$. 
However, with our restricted assumption (2.1.2) - (2.1.4) 
we cannot hope for such a nice form. 
\par Near $(0,\1,0)$ the shape of $Q$ is essentially determined 
by the linear involution $\sigma\:= d\tau(0,\1, 0)$
on  $E={\Cal T}_{(0,\1,0)}{\Cal U}\simeq {\Cal T}_mM$. 

\par We will use  the natural identification 
$E=\s_m^*\times \s_m\times V$. 
Define $W=\s_\m^* \times \s_m \times V_{\fix}$. Then 
it follows from (2.2.2) that 
$E=W\oplus V_{\eff}$ is the decomposition 
of $E$ into fixed and effective part 
of the isotropy action of $T_m$ on $E\simeq {\Cal T}_m M$.
Then (2.1.2) implies that the 
involution $\sigma$ leaves the decomposition 
$E=W\oplus V_{\eff}$ invariant. 
Hence 
$$\sigma=\pmatrix \sigma\res_W & 0 \\ 0 & \sigma\res_{V_{\eff}}
\endpmatrix.\leqno(2.2.4)$$
Accordingly we have a splitting 
$${\Cal T}_{(0,\1,0)} Q=\left({\Cal T}_{(0,\1,0)} Q\cap W\right)
\oplus \left({\Cal T}_{(0,\1,0)} Q\cap V_{\eff}\right).$$
Next we will analyze the pieces $\sigma\res_W$ and 
$\sigma\res_{V_{\eff}}$. We start with $\sigma\res_{V_{\eff}}$. 
Notice that it follows from (2.1.2) that 
$$\pi(t)\circ \sigma\res_{V_{\eff}}= \sigma\res_{V_{\eff}}
\circ \pi(t^{-1})\qquad (t\in T_m) . \leqno (2.2.5)$$
Thus we can apply Lemma 2.1.2 (to $V=V_{\eff}$, 
$\tau=\sigma\res_{V_{\eff}}$ and $V_1={\Cal T}_{(0,\1,0)} Q\cap 
V_{\eff}$) and conclude:
$$\sigma\res_{V_{\eff}}\quad \hbox{ is anti-symplectic}. \leqno 
(2.2.6)$$
For the scalar product $\la\cdot, \cdot\ra$ on $V_{\eff}$
this means that we may assume in addition 
that it is invariant under $\sigma
\res_{V_{\eff}}$. 
\par Write $V_{\eff}=V_{\eff, 1}\oplus V_{\eff, -1}$
for the decomposition of $V$ into eigenspaces 
of $d\sigma\res_{V_{\eff}}$. Thus 

$$d\sigma\res_{V_{\eff}}=\pmatrix \id_{V_{\eff, 1}} & 0\\
0 & -\id_{V_{\eff, -1}}\endpmatrix.\leqno(2.2.7)$$
Notice that $V_{\eff,1}={\Cal T}_{(0,\1,0)} Q\cap V_{\eff}$. 

\ssk Next we turn our attention to 
$\sigma\res_W$. From (2.1.5) and the concrete formula
(2.2.2) for the $S_m$-action  it follows that 
$$\sigma\res_W= \pmatrix * & 0 & *\\ * &-\id_{\s_m} &*\\
* &0&*\endpmatrix \leqno(2.2.8)$$
with respect to a basis of $W$ compatible 
with  $W=\s_m^*\times 
\s_m
\times V_{\fix}$.  
A consequence of (2.1.3) is 
$d\Phi(0,\1,0)\circ \sigma=d\Phi(0,\1,0)$. 
From (2.2.4) it hence follows that 
$d\Phi(0,\1,0)\res_W\circ \sigma_W=d\Phi(0,\1,0)\res_W$. 
As $d\Phi(0,\1,0)\res_W$ is the projection 
$W\to \s_m^*$ along $\s_m\times V_{\fix}$, the pattern
(2.2.8) further simplifies to 
 
$$\sigma\res_W= \pmatrix \id_{\s_m^*} & 0 & 0\\ * &-\id_{\s_m} &*\\
* &0&D\endpmatrix \leqno(2.2.9)$$
for some linear operator $D\: V_{\fix}\to V_{\fix}$. 
As $\sigma\res_W\circ \sigma\res_W=\id_W$ we derive 
from (2.2.9) that $D^2=\id_{V_{\fix}}$. Accordingly
we obtain an eigenspace decomposition 
$V_{\fix}=V_{\fix,1}\oplus V_{\fix, -1}$ for  $D$. 
With respect to the refined decomposition 
$W=\s_m^*\times \s_m \times V_{\fix,1}\times V_{\fix, -1}$
we then have 
$$\sigma\res_W= \pmatrix \id_{\s_m^*} & 0 & 0 &0\\ A &-\id_{\s_m} &
B_1& B_2\\
C_1 & 0  &\id_{V_{\fix, 1}}& 0\\ 
C_2 & 0 & 0 & -\id_{V_{\fix,-1}}\endpmatrix.$$
Employing again $\sigma\res_W\circ \sigma\res_W=\id_W$
we obtain that  $B_2=0$ and $C_1=0$. Thus 
$$\sigma\res_W= \pmatrix \id_{\s_m^*} & 0 & 0 &0\\ A &-\id_{\s_m} &
B& 0\\
0 & 0  &\id_{V_{\fix, 1}}& 0\\ 
C& 0 & 0 & -\id_{V_{\fix,-1}}\endpmatrix\leqno(2.2.10)$$
for linear operators $A\: \s_\m^*\to \s_m$, 
$B\: V_{\fix,1}\to \s_m$ and $C\: \s_m^*\to 
V_{\fix, -1}$. 
Combining (2.2.4) and (2.2.10) we then obtain 
$$\sigma= 
\pmatrix \id_{\s_m^*} & 0 & 0 &0 & & \\ A &-\id_{\s_m} &
B& 0 & & \\
0 & 0  &\id_{V_{\fix, 1}}& 0 & & \\ 
C& 0 & 0 & -\id_{V_{\fix,-1}}  & & \\ 
& & & & \id_{V_{\eff,1}} & 0\\ 
& & & & 0  & -\id_{V_{\eff, -1}}\endpmatrix\leqno(2.2.11)$$
with respect to a basis of $E$ compatible with 
$E=\s_m^*\times \s_m \times V_{\fix, 1}\times V_{\fix, -1}
\times V_{\eff, 1}\times V_{\eff, -1}$. 
\par With the help of (2.2.11) we can now
determine the tangent space ${\Cal T}_{(0,\1,0)}Q\subeq E$. 
Notice that 
$${\Cal T}_{(0,\1,0)} Q=\{ v\in E\: \sigma(v)=v\}$$
so that (2.2.11) implies that 
$${\Cal T}_{(0,\1,0)} Q=\left\{ \pmatrix x \\ {1\over 2} 
(Ax+By)\\  y\\ {1\over 2}Cx\\ z\\ 0\endpmatrix\: 
x\in \s_m^*, \ y\in V_{\fix, 1}, z\in V_{\eff, 1}\right\}.
\leqno(2.2.12)$$
Write $\pr\: E\to \s_m^*\times V_{\fix,1}\times V_{\eff,1}$
for the projection along 
$\s_m\times V_{\fix,-1}\times V_{\eff,-1}$. Then (2.2.12)
implies that $\pr\res_{{\Cal T}_{(0,\1,0)} Q}\: 
{\Cal T}_{(0,\1,0)}\to \s_m^*\times V_{\fix,1}\times V_{\eff,1}$ 
is a linear isomorphism. 
This in turn allows us to apply the implicit function theorem: 
there exist an open neigborhood $U_1$ of $0$ in 
$\s_m^*\times V_{\fix,1}\times V_{\eff,1}$, an open 
neigborhood $U_2$ of $(\1,0,0)$ in $S_m\times V_{\fix,-1}\times 
V_{\eff, -1}$ and a differentiable map 
$$\psi=(\psi_S, \psi_{\fix},\psi_{\eff})\: U_1\to U_2$$
such that $\psi(0,0,0)=(\1,0,0)$ and 
$$Q\cap (U_1\times U_2)=\left\{ 
\pmatrix x\\ \psi_S(x, y, z)\\  y\\  \psi_{\fix}(x,y,z)\\  
z\\ \psi_{\eff}(x,y,z))\endpmatrix\in {\Cal T}^* S_m \times V\: 
(x,y,z)\in U_1\right\}.\leqno(2.2.13)$$
One can say a little bit more about the map 
$\psi$ when we notice that ${\Cal T}_{(0,\1,0)} Q$ can 
equally be expressed as 
$${\Cal T}_{(0,1,0)}Q=\left\{\pmatrix u \\ d\psi_S(0)(u,v,w)\\
v \\ d\psi_{\fix}(0)(u,v,w) \\ w\\ d\psi_{\eff}(0)(u,v,w)
\endpmatrix\in E\: u\in \s_m^*, \ v\in V_{\fix, 1}
,\ w\in V_{\eff, 1}\right\}.\leqno(2.2.14)$$
Comparing (2.2.12) with (2.2.14) yields
$$d\psi_{\eff}(0)=0.\leqno(2.2.15)$$

\par We are now ready to summarize the discussion 
of this subsection. In more compact notation 
we have proved the following: 

\proclaim{Theorem 2.2.2} {\rm (Local normal forms
for $\Phi$ and $\Phi\res_Q$)} Let $(M,\omega)$
be a connected symplectic manifold endowed 
with a Hamiltonian torus action and an involution 
$\tau\: M\to M$ which satisfy 
{\rm (2.1.2) - (2.1.4)}. Let $m\in Q$. Write 
$k=\dim \s_m=\dim \t/\t_m$ and identify 
$\t^*\simeq \R^k\times \t_m^*$. 
Then there exist an open neighborhood $U$ of $m$ and 
symplectic coordinates $\x,\y, {\text {\bf q}}, 
{\text {\bf p}} = x_1, \ldots, x_k, y_1, \ldots, y_k, 
q_1, \ldots, 
q_N, p_1, \ldots, p_N$ with 
$\x(m)=\y(m)={\text {\bf q}}(m)= {\text {\bf p}}(m)=0$
such that:
\roster
\item On $U$ the momentum map $\Phi\:M \to 
\t^*\simeq \R^k\times \t_m^*$ is
given by 
$$\Phi=\Phi(m)+ \Big(\x, {1\over 2}\sum_{j=1}^l \lambda_j(q_j^2
+p_j^2)\Big)$$ 
where $l\leq N$ and $\lambda_1, \ldots,\lambda_l\in 
\t^*\bs \{0\}$. 
\item The restriction  
$(\x, {\text {\bf q}})\: 
Q\cap U\to \R^n$ of $(\x,{\text {\bf q}})$
to $Q\cap U$ is a diffeomorphism 
onto an open ball $B_r^n(0)$ of radius $r>0$ in $\R^n$. Furthermore 
$$\Phi\res_Q= \Phi(m)+ \Big(\x, 
{1\over 2}\sum_{j=1}^l 
\lambda_j(q_j^2 +\psi_j(\x, {\text {\bf q}})^2)\Big) ,
$$
where $\psi=(\psi_1, \ldots, \psi_l)\: B_r^n(0)\to \R^l$
is a differentiable map 
with $\psi(0)=0$ and $d\psi(0)=0$. 
\endroster 
\endproclaim 
\demo{Proof} To explain the notation in the theorem:
$(\x, \y)$ are symplectic coordinates for 
${\Cal T}^* S_m =\s_m^*\times S_m$ with 
$\x$ corresponding to $\s_m^*$ and $\y$ to $S_m$; 
next $({\text{\bf q}}, {\text{\bf p}})$ are symplectic 
coordinates for $V$ compatible 
with the weight space decomposition $V=\bigoplus_{\lambda}
V_\lambda$ and moreover $q_1, \ldots, q_l$ corresponding 
to $V_{\eff, 1}$ and $q_{l+1}, \ldots, q_N$ corresponding 
to $V_{\fix, 1}$ (and similar for  ${\text{\bf p}}$). 
\par The expression for $\Phi$ in 1. then follows from 
(2.2.3) and (2.2.6) which implied that 
the inner product on $V_{\eff}$ could be chosen 
$\sigma\res_{V_{\eff}}$-invariant. Finally, the 
assertion in 2. follows from 1. and (2.2.13) combined 
with (2.2.15). Here the map $\psi$ corresponds to 
$\psi_{\eff}$ above.\qed 
\enddemo

\head 3. The convexity theorem \endhead

The objective of this subsection is to prove 

\proclaim {Theorem 3.1} Let $M$ be a compact connected 
symplectic manifold with Hamiltonian torus action 
$T\times M\to M$ 
and momentum map $\Phi\: M\to \t^*$. 
Further let $\tau\: M\to M$ 
be an involutive diffeomorphism with 
fixed point set $Q$ such that 
\roster
\item $t\circ \tau=\tau\circ t^{-1}$ for all $t\in T$. 
\item $\Phi\circ \tau=\Phi$. 
\item $Q$ is a Lagrangian submanifold of $M$. 
\endroster 
Then $\Phi(Q)=\Phi(M)$. In particular 
$\Phi(Q)$  is a convex subset of $\t^*$.
Moreover, the same assertions hold if $Q$ is replaced with
any of its connected components.
\endproclaim

Our arguments follow the approach of Duistermaat (cf. [4, Sect. 2]). 
However, they have to be adapted to the more general case where 
$\tau$ is not necessarily anti-symplectic.

\par From now on we will assume that 
$M$ is compact. In particular, $Q$ is a compact submanifold of $M$.  
Furthermore, we will require that $T$ acts freely, i.e., 
$\bigcap_{m\in M} T_m=\{\1\}$. But 
notice that this is not a severe restriction 
as we can always replace 
$T$ with $T/ \bigcap_{m\in M} T_m$. 

\par The key result toward convexity of $\Phi(Q)$ is the 
following central fact which 
generalizes [6, Lemma 32.1] and [4, Lemma 2.4].

\proclaim {Proposition 3.2} Let $X\in \t$. Then the function 
$\Phi_X\res_Q$ has a unique local maximal value. 
\endproclaim 

\demo{Proof} Fix $X\in \t$ and let $m\in Q$ be a critical point of $\Phi_X\res_Q$. 
By Lemma 2.2.1 we know that $m$ is fixed by $\exp(\R X)$. Thus 
replacing $T$ by $\oline{\exp(\R X)}$ we may assume that $d(\Phi\res_Q)(m)=0$. 
In particular we obtain from Theorem 2.2.2  with 
$\s_m	=\{0\} $ that 
$$\Phi_X=\Phi_X(m) +{1\over 2}\sum_{j=1}^l 
\lambda_j(X) (q_j^2 +p_j^2)\leqno(3.1)$$
and 
$$\Phi_X\res_Q= \Phi_X(m) +{1\over 2}\sum_{j=1}^l 
\lambda_j(X) (q_j^2 +\psi_j({\text {\bf q}})^2)\leqno(3.2)$$
hold in a neighborhood of $U$ of $m$. 
\ssk\nin 
{\it Claim:}   $m$ is a local maximum for $\Phi\res_Q$
iff $\lambda_j(X)\leq 0$ for all $1\leq j\leq l$. 
\ssk 
In view of (3.1) this is clear if all $\lambda_j(X)\leq 0$. 
To obtain the other direction assume that $\lambda_j(X)>0$ for some $j$. 
W.l.o.g. we may assume that $\lambda_1(X)>0$. 
For  each $j$ we then have 
$$\psi_j(q_1, 0, \ldots, 0)= q_1^2 h_j(q_1)$$
for a continuous function $h_j$. This is because of 
$\psi_j(0)=d\psi_j(0)=0$
for all $1\leq j\leq l$ (see Theorem 2.2.2 (2)). 
Thus for small and non-zero $q_1$ we obtain from 
(3.2)

$$\eqalign{\Phi_X\res_Q(q_1, 0,\ldots, 0) 
&= \Phi_X(m) +{1\over 2}\Big(\lambda_1(X) q_1^2 
(1 +q_1^2 h_1(q_1)^2) \cr 
&  + q_1^4 \sum_{j=2}^l \lambda_j(X) h_j(q_1)^2\Big)\leq \Phi_X(m).\cr}$$ 
Dividing by ${1\over 2}q_1^2$ implies that 

$$\lambda_1(X) (1 +q_1^2 h_1(q_1)^2) 
+ q_1^2 \sum_{j=2}^l \lambda_j(X) h_j(q_1)^2 \leq 0$$
for $q_1$ small and non zero. 
This clearly contradicts $\lambda_1(X)>0$ and proves our claim. 
\par It follows from our claim and (3.1) and (3.2) 
that $m$ is a local maximum 
for $\Phi_X\res_Q$ if and only if $m$ is a local maximum for 
$\Phi_X$. As $M$ is compact, 
Morse theoretic arguments imply that the set of points of $M$ where $\Phi_X$ attains
a local maximum is connected (see the proof of 
[6, Cor. 32.1] or [4, Lemma 2.4]).
This completes the proof of the proposition.
\qed\enddemo

\proclaim {Lemma 3.3} If $m\in Q$ is such that $\xi=\Phi(m)$ is 
a boundary point of $\Phi(Q)$, then $\t_m\neq \{0\}$. \endproclaim 

\demo{Proof} If $\t_m=\{0\}$, then (2.2.1) implies that $d(\Phi\res_Q)(m)$ is 
surjective, contradicting the fact that $\xi$ is a boundary point. \qed 
\enddemo

\par Let $\R^+=]0,\infty[$ and $\R_0^+=[0,\infty[$ and 
define a closed  convex cone in $\t_m^*$ by  
$$ \Gamma_m=\{ \sum_{j=1}^l s_j\lambda_j\: s=(s_1, \ldots, s_l)\in (\R_0^+)^l\}.$$

\proclaim{Lemma 3.4}
If $\Gamma_m=\t_m^*$, then $\im \Phi \res_Q$ contains an open
neighborhood of $\Phi(m)$ in $\t_m^*$.
\endproclaim

\demo{Proof}
Since $\Gamma_m=\t_m^*$, there exists a vector
$$\v=(0;v_1,\dots,v_l) \in \R^k \times (\R^+)^l \quad \hbox{such that}
\quad \sum_{j=1}^n v_j \lambda_j=0\in \t_m^* . \leqno(3.3)
$$
Then for any neighborhood $U$ of $\v\in \R^k \times (\R^+)^l $ the set
$$\{ (u_1,\dots u_k,u_{k+1}\lambda_1 + \dots + u_n \lambda_l) :
(u_1,\dots,u_n)\in U \} $$
contains an open neighborhood of $0$ in $\t^*$.
\par Recall the map $\psi:B_r^n(0) \to \R^l $ from Theorem
2.2.2(2). We define 
$$\Psi\: B_r^n(0) \to \R^k \times (\R_0^+)^l ,\ \ (\x,\q) \mapsto
(\x,{1\over 2} (q_1^2 +\psi_1(\x,\q)^2), \ldots, {1\over 2}(q_l^2 +
\psi_l(\x,\q)^2)) . \leqno(3.4) $$
According to Theorem 2.2.2(2) it is sufficient to show that
$\im \Psi$ contains a point $\v$ as in (3.3) as inner point.
Fix $\v$ satisfying (3.3). The theorem will be proved if we can show 
that $\im \Psi$ contains
$s\v$ as an inner point for some $s>0$. In the following we will verify 
this claim. We may assume that $\|\v\|=1$.

\par The properties $\psi(0)=0, \ d\psi(0)=0 $ from Theorem 2.2.2(2)
are crucial. Together with Taylor's formula they imply that we can find
a constant $K>0$ such that for all $1\leq j\leq l$,
$$ |\psi_j(\x,\q)| \leq K \|(\x,\q)\| ^2 , \leqno(3.5) $$
for every $(\x,\q)\in \R^k \times (\R^+)^l $ sufficiently close to 
$(0,0)$.

\par There is another constant $C>0$ such that the ball $B_{sC}^n(s\v)$
lies entirely in $\R^k \times (\R^+)^l $ for all $s>0$. We notice that 
every point in
$B_{sC}^n(s\v)$ can be written as $(x_1,\dots,x_k,{1\over 2} q_1^2,\dots
,{1\over 2} q_l^2)$ with unique $x_1,\dots,x_k \in \R , q_1,\dots,q_l
\in \R^+ $. From now on let $0<s\leq 1$. Observe that  
the condition 
$$(x_1,\dots,x_k,{1\over 2} q_1^2,\dots
,{1\over 2} q_l^2) \in B_{sC}^n(s\v)$$ 
puts restrictions on the vector
$(x_1,\dots,x_k,q_1,\dots,q_l)$: there is a constant $\tilde C >0$ (independent 
of the vector)  such that 
$$\| (x_1,\dots,x_k,q_1,\dots,q_l) \| \leq \sqrt{s}\tilde{C} . $$
In particular, (3.5) implies that for some $\tilde{K}>0$,
$$ \| (0,\dots,0,{1\over 2} \psi_1(\x,\q)^2,\dots,{1\over 2} 
\psi_l(\x,\q)^2 \| \leq \tilde{K} s^2 \leqno(3.6) $$
for all $(\x,\q)=(x_1,\dots,x_k,q_1,\dots,q_l)$ with
$(x_1,\dots,x_k,{1\over 2} q_1^2,\dots,{1\over 2} q_l^2) \in 
B_{sC}^n(s\v)$.

\par Choose  $0<s_0\leq 1$ small enough  such that $\tilde{K}s_0^2 < s_0 C$ holds. 
Set $\epsilon=s_0 C$. We are now in a position to apply Brouwer's fixed point theorem:
consider the mapping
$$ \Xi\: \oline{B_{\epsilon}^n(0)} \to \oline{B_{\epsilon}^n(0)}, $$
$$(x_1,\dots,x_k,{1\over 2} q_1^2,\dots,{1\over 2} q_l^2)-s_0 \v 
\mapsto -(0,\dots,0,{1\over 2} \psi_1(\x,\q)^2,\dots,{1\over 2} 
\psi_l(\x,\q)^2) . $$
If $(\x,{1\over 2} q_1^2,\dots,{1\over 2} q_l^2)-s_0\v $ is a
fixed point of $\Xi$, then
$$ \Psi(\x,\q)=(\x,{1\over 2} (q_1^2 +\psi_1(\x,\q)^2), \ldots, {1\over 2}(q_l^2 +
\psi_l(\x,\q)^2))=s_0 \v . $$

\par We want to show that $s_0\v$ is an inner point of $\im \Psi$.
Notice that by choosing  $s_0$ small enough we can assume that the 
point $(\x,\q)$ with $\Psi(\x,\q)=s_0\v $ is arbitrarily close to
$(0,0)$.
The mapping $\Psi$ is submersive at most points close enough to the origin
as a look at its derivative shows:
$$ \eqalign{d\Psi(\x,\q)&=\pmatrix 
\id_k & & & \\ & q_1 & & \\ & &\ddots& \\ & & & q_l\endpmatrix
\cr 
\qquad & +\pmatrix 0 & 0 \\
(\psi_i(\x,\q) {\partial\psi_i\over \partial x_j}(\x,\q) )_{i,j} &
(\psi_i(\x,\q) {\partial\psi_i\over \partial q_j}(\x,\q) )_{i,j}
\endpmatrix\cr}$$
Relation (3.5) implies that for $(\x,\q)$ approaching $(0,0)$ the 
entries in the second summand become arbitrarily small compared 
to those in the first summand. Since the $\q$ under consideration satisfy 
$q_1,\dots,q_l >0 $, we conclude
that 
$$ \det(\Psi(\x,\q))>0 \quad \hbox{for} \ (\x,\q) \ \hbox{close enough to} \ (0,0) . $$
This finishes the proof.
\qed
\enddemo

If $C$ is a closed convex cone in an Euclidean  vector space $E$, then its 
{\it dual cone} in $E^*$ is defined by 
$$C^\star=\{ \lambda\in E^*\: \lambda\res_C\geq 0\}.$$
Recall that $C^\star$ is a closed convex cone in $E^*$. One has 
$(C^\star)^\star=C$, and in particular $C^\star=\{0\}$ if and only if 
$C=E$.

\proclaim {Lemma 3.5} Let $\xi\in \Phi(Q)$ be a boundary 
point and $m\in Q$ such that $\Phi(m)=\xi$. 
Then the following assertions hold:
\roster
\item $\Gamma_m\neq \t_m^*$. In particular $\{0\}\neq \Gamma_m^\star\subeq \t_m$. 
\item For all $0\neq X\in -\Gamma_m^\star\subeq \t_m$, 
one has $\Phi_X\res_Q\leq \Phi_X(m)=\la
\xi , X\ra$. In particular $\im \Phi$ is 
contained in the half space $\{ \lambda\in \t^*\:\ (\xi-\lambda)(X)\geq 0\}$. 
\item $\xi$ is a boundary point of $\Phi(M)$.
\endroster
\endproclaim

\demo{Proof} 1. Suppose that $\Gamma_m=\t_m^*$. Then Lemma 3.4
implies that $\im \Phi$ contains a neighborhood of $\xi$.
 But this contradicts the fact that $\xi$ is a boundary point. 
\par\nin 2. According to 1. we have $\Gamma_m^\star\neq \{0\}$. Let $0\neq X\in 
-\Gamma_m^\star$. Then (3.2) implies that 
$\Phi_X(m')\leq \Phi_X(m)=\xi(X)$
for all $m'$ in a small neighborhood of $m$ in $Q$. Thus 
$\Phi_X(m)$ is a local maximal value of $\Phi_X\res_Q$. Hence Proposition 
3.2 implies that $\Phi_X\res_Q\leq \Phi_X(m)$. This completes 
the proof of 2. 
\par\nin 3. A slight modification (refer to (3.1) instead of (3.2))
of the argument just given shows that
there is a $0\neq X\in -\Gamma_m^\star$ such that
$$ \Phi_X(m') \leq \Phi_X(m)=\la \xi , X\ra \qquad \hbox{for all $m'$ in
a neighborhood of $m$ in $M$.} $$
The function $\Phi_X$ on $M$ has a unique local maximal value (see [4,
Lemma 2.4]). Therefore, $\Phi(M)$ must lie entirely in the halfspace
$\{ \lambda\in \t^*\:\ (\xi-\lambda)(X)\geq 0\}$, implying that $\xi$
is a boundary point of $\Phi(M)$.
\qed 
\enddemo

Let us define {\it regular elements in $Q$} by 
$$Q_{\reg}=\{ m\in Q\: d\Phi\res_Q(m) \  \hbox{is 
surjective}\}.$$

\proclaim{Lemma 3.6} The following assertions hold: 
\roster
\item $Q_{\reg}$ is open and dense in $Q$. 
\item The interior $\Int \Phi(Q)$  is 
dense  in $\Phi(Q)$. 
\endroster
\endproclaim 
\demo{Proof} 1. It is clear that $Q_{\reg}$ is open. Let us 
show that $Q_{\reg}$ is dense. For that fix $m\in Q\bs Q_{\reg}$. 
Then $d\Phi\res_Q(m)$ is not onto and hence $\t_m\neq \{0\}$
by (2.2.1). In the $(x,\q)$-coordinates near $m$ we have 
$$\Phi_Q(\x,\q)=\Phi(m)+ \Big(\x, {1\over 2} 
\sum_{j=1}^l (q_j^2 +\psi_j(\x,\q)^2)\lambda_j\Big)$$
(cf.\ Theorem 2.2.2(2)). 
As we assume that $\bigcap_{m\in M}T_m=\{\1\}$ it follows that 
$\lambda_1, \ldots, \lambda_l$ linearly span $\t_m^*$.
Recall the definition of the map $\Psi$ from (3.4).
As was shown at the end of the proof of Lemma 3.4
we can find elements $(0,\q)$ arbitrarily close 
to $m=(0,0)$ such that $d\Psi\res_Q(0,\q)$ is invertible. This 
in turns implies that $d\Phi\res_Q(0,\q)$ is surjective and completes 
the proof of 1. 
\par\nin 2. This is immediate from 1. \qed 
\enddemo

\msk\nin{\bf Proof of Theorem 3.1.} We first show 
that $\Phi(Q)$ is convex. For that notice that 
$\Phi(Q)$ is a compact subset 
of $\t^*$ with open and dense interior 
(Lemma 3.6(2)).  
By [12, Satz 3.3]
$\Phi(Q)$ will be convex if each boundary point 
lies on a half space 
containing $\Phi(Q)$. In view of Lemma 3.5(2)  this is satisfied and so 
$\Phi(Q)$ is convex.

\par Next we show that $\Phi(Q)=\Phi(M)$. 
Assume that some extremal point $\eta$ of 
$\Phi(M)$ does not 
lie in $\Phi(Q)$. From Lemma 3.6 we know there is a point 
$\zeta$ in
$\Int \Phi(Q)$. On the line segment connecting 
$\eta$ and $\zeta$ there 
must be a boundary point $\beta$ of $\Phi(Q)$. Since 
$\Phi(M)$ is a convex 
polyhedron (and since $\eta \not\in \Phi(Q)$), we see that 
$\beta$ must be contained
in $\Int \Phi(M)$. This contradicts Lemma 3.5(3).
\qed

\head 4. Application to Kostant's Theorem \endhead

In this section we will use Theorem 3.1 to give a
symplectic proof of Kostant's non-linear convexity 
theorem. We start with the introduction of the 
necessary notation and the statement of Kostant's
result. 

\subhead 4.1 Notation \endsubhead 

Let $G$ denote a connected semisimple Lie group. Universal 
complexifications of Lie groups will denoted by a subscript $\C$, 
i.e. $G_\C$ is the universal complexification of $G$ etc.  
For what follows it is no loss of generality when we 
assume that $G\subeq G_\C$ and that $G_\C$ is simply connected. 

\par Write $\g$ for the Lie algebra of $G$. Complexifications 
of Lie algebras shall be denoted by the subsript $\C$, i.e.
$\g_\C$ is the complexification of $\g$ etc. 

\par Let $\g=\k+\p$ be a Cartan decomposition of $\g$ with 
$\k$ a maximal compact subalgebra. Fix a maximal 
abelian subspace $\a\subeq \p$  and denote by $\Sigma=\Sigma(\g,\a)$
the corresponding restricted root system in $\a^*$, the dual 
of $\a$. For each $\alpha\in \Sigma$ let $\g^\alpha=\{ Y\in \g\: 
(\forall H\in \a)\ [H,Y]=\alpha(H)\}$ the associated 
root space. With $\m=\z_\k(\a)$ one then has the root space 
decomposition 
$$\g=\a+\m +\bigoplus_{\alpha\in \Sigma} \g^\alpha.$$  
Select a positive system $\Sigma^+\subeq \Sigma$ 
and define the nilpotent subalgebra $\n=\bigoplus_{\alpha\in 
\Sigma^+}\g^\alpha$. 

\par On the group level we denote by $A, K$ and $N$
the analytic subgroups of $G$ with Lie algebras 
$\a, \k $ and $\n$. Then there is 
the Iwasawa decomposition of $G$ which states that 
the multiplication mapping 
$$N\times A\times K\to G, \ \ (n, a, k)\mapsto nak$$
is an analytic diffeomorphism. For $g\in G$ let us 
denote by $\tilde a(g)$ the $A$-component of $g$ in the Iwasawa 
decomposition. 

\par Set $M=Z_K(\a)$ and note that the Lie algebra 
of $M$ is $\m$. The Weyl group of $\Sigma$ can then be defined 
by ${\Cal W}=N_K(\a)/ M$. 

\par We are ready to state Kostant's theorem [11]: 

\proclaim {Theorem 4.1.1} Let $Y\in \a$.  Then 
$$\log \tilde a(K\exp(Y))=\conv ({\Cal W}.Y)$$
where $\conv(\cdot)$ denotes the convex hull
of $(\cdot)$. \qed \endproclaim

\subhead 4.2 Symplectic methods for the complex case\endsubhead

In case $G$ is complex a symplectic proof of Theorem 4.1.1 was 
given by Lu and Ratiu [13]. The objective 
of this section is to briefly recall their method. 

\par Let us assume that $G$ is complex, i.e. $\g$ is a complex Lie algebra. 
Then the Cartan decomposition of $\g$ is given by 
$\g=\k+i\k$, i.e.  $\p=i\k$. Furthermore $\a=i\t$ with 
$\t$ a maximal toral subalgebra  in $\k$. Set $T=\exp(\t)$. 

\par Define a solvable subalgebra of $\g$ by $\b=\a+\n$ (despite 
the notation this is not a Borel subalgebra of $\g$). Write $B=AN$ 
for the corresponding group and notice that $B$ is invariant under 
conjugation by the torus $T$. 

\par Let us denote by $\kappa$ the Cartan-Killing form 
of the complex Lie algebra $\g$ and define a symmetric 
$\R$-valued bilinear form on $\g$ 
by 
$$ \B\: \g\times \g\to \R, \ \ \B(X, Y)=\Im \kappa(X, Y).$$
Notice that $\B$ is invariant and non-degenerate. The 
important fact is that both $\b$ and $\k$ are isotropic 
for $\B$; in other words $(\g, \b, \k)$ becomes a 
Manin-triple [3]. 
Likewise $(G, B, K)$ is a Manin-triple. Recall that this 
implies that $B\simeq G/K$ carries a natural structure 
of a Poisson Lie group [3]. In order to describe the symplectic 
leaves write $\tilde b(g)$ for the $B$-part of $g\in G$ in the decomposition 
$G=B\cdot K$. Then the symplectic leaves in $B$ are given by 
$M_a=\tilde b(Ka)$ for $a\in A$.  
As manifolds $M_a\simeq K/K_a$ with 
$K_a=Z_K(a)$. Notice that $K$ does not act symplectically on $M_a$; 
however $T$ does and the $T$-action is Hamiltonian. It was 
established in [13] that the corresponding momentum map 
is the non-linear Iwasawa projection:
$$\Phi\: M_a\to \a\simeq \t^*, \ \ \tilde b(ka)\mapsto \log \tilde a(ka).
\leqno(4.2.1)$$  
Standard structure theory implies that $\Fix (M_a)={\Cal W}.a$. 
Thus (4.2.1) combined with the Atiyah-Guillemin-Sternberg convexity 
theorem gives a symplectic proof of  Theorem 4.1.1 
in the case of $G$ complex [13].

\par For later reference we give an explicit formula
for the symplectic form on a leaf $M_a=\tilde b(Ka)$. 
For $X\in \k$ let us denote by $\tilde X$ the corresponding 
vector field on $M_a$, i.e.
$$\tilde X_b={d\over dt}\Big|_{t=0} \tilde b(\exp(tX)b )\qquad 
(b\in M_a).$$
Then the symplectic form $\omega$ of $M_a$ is given by 

$$ \omega_b(\tilde X_b, \tilde Y_b)=\B (\pr_\k(\Ad(b)^{-1}X), 
\Ad(b)^{-1}	 Y) \qquad (b\in M_a; \ X, Y\in \k)\leqno (4.2.2)$$
where $\pr_\k\: \g \to \k$ is the projection along $\b$. 
This is immediate from [3, (11.1.2)].

\subhead 4.3 Symplectic methods for the real case\endsubhead 

Define a maximal compact subalgebra 
$\uu$ in $\g_\C$ by $\uu=\k+i\p$ and let $U$ be the 
corresponding maximal compact subgroup of $G_\C$.
Notice that  $\g_\C=\uu+i\uu$ is a Cartan decomposition of $\g_\C$. 
If $\t_1$ denotes a maximal torus in $\m=\z_\k(\a)$, then 
$\a_1=\a+i\t_1$ defines a maximal abelian subspace of $i\uu$. 
Write $\Sigma_1=\Sigma_1(\g_\C, \a_1)$ for the corresponding 
root system. Fix a positive system $\Sigma_1^+$ of 
$\Sigma_1$. Without loss of generality we may assume 
that $\Sigma_1^+$ and $\Sigma^+$ are compatible, i.e.
$\Sigma_1^+\res_\a\subeq \Sigma^+\cup\{0\}$.  
Write $\n_1=\bigoplus_{\alpha\in \Sigma_1^+} \g_\C^\alpha$
for the nilpotent subalgebra of $\g_\C$ associated 
to $\Sigma_1^+$. Likewise we denote by $N_1$ the corresponding 
subgroup of $G_\C$. Notice that $N_\C\subeq N_1$ but 
generally $N_\C\subsetneq N_1$ unless $\m=\t_1$ is abelian. 
Clearly we have $A\subeq A_1$ and so $B\subeq B_1$ where 
$B=AN$ and $B_1=A_1N_1$. Finally let us define the torus
$T=\exp(i\a)$. 

\par Fix $a\in A\subeq A_1$ and consider the symplectic 
manifold $M_a=\tilde b(Ua)$ of the Poisson Lie group $B_1$ 
(cf.\ Subsection 4.2). As explained in Subsection 4.2
the action of $T$ on $M_a$ is Hamiltonian with momentum 
map $\Phi\: M_a\to \a\simeq \t^*$
given by (4.2.2)
$$\Phi(\tilde b_1(ua))=\log \tilde a(ug)\qquad (u\in U)$$
with $\tilde a(g)$ the $A$-part of $g\in G_\C$ in the 
decomposition $G_\C=N_1 A\exp(i\t_1) U$. 
\par Set $Q_a=\tilde b(Ka)$. 
Then the restriction of $\Phi$ to $Q_a$ is the nonlinear
Iwasawa projection, i.e $\Phi(\tilde b(ka))=\log \tilde a(ka)$ for 
all $k\in K$. We are interested in the image $\Phi(Q_a)$. 

\par Denote by $\tau\: G_\C\to G_\C$
the complex conjugation with respect to the real form $G$. 
Notice that $\tau$ is an involutive diffeomorphism of $G_\C$ 
which induces  an involution on $M_a$, say $\tau_a$, by the 
prescription 
$$\tau_a(\tilde b(ua))=\tilde b(\tau(ua))=\tilde b(\tau(u)a) \qquad (u\in U).$$
Standard stucture theory shows that the connected 
component containing $a\in M_a$ of the fixed-point set of the involution 
$\tau_a\: M_a\to M_a$ is given by $Q_a=\tilde b(Ka)$. 
We collect some important properties of $\tau_a$ and $Q_a$:

\proclaim {Lemma 4.3.1} The following assertions hold: 
\roster 
\item $Q_a\subeq M_a$ is Lagrangian submanifold of $M_a$. 
\item For all $t\in T$ and $m\in M_a$ one has 
$\tau_a(t.m)=t^{-1}.\tau_a(m)$. 
\item $\Phi\circ \tau_a=\Phi$. 

\endroster 
\endproclaim 
 
\demo{Proof} 1. It is sufficient to show that 
$\omega_b(\tilde X_b, \tilde Y_b)=0$ for all 
$b\in Q_a$ and $X, Y\in \k$. Notice that  $B$ and $\k$ are 
$\tau$-fixed. Thus (4.2.2) implies
$$\eqalign{\omega_b(\tilde X_b, \tilde Y_b) & =\B\left(\pr_\uu(\Ad(b)^{-1}X), 
\Ad(b)^{-1}Y\right)\cr 
&=\B \left(\pr_\k (\tau(\Ad(b)^{-1}X)), 
\tau(\Ad(b)^{-1}Y)\right)\cr
&=\B \left(\tau(\pr_\k (\Ad(b)^{-1}X)), 
\tau(\Ad(b)^{-1}Y)\right).\cr}$$
By the definition of $\B$ we have $\B\circ (\tau\times\tau)=-\B$. Thus 
our computation above gives 
$\omega_b(\tilde X_b, \tilde Y_b)=-\omega_b(\tilde X_b, \tilde Y_b)$; 
completing the proof of 1. 
\par\nin 2. Let $b\in M_a$. Then $b=\tilde b(ua)$ for some 
$u\in U$. Hence for $t\in T$: 
$$\tau_a(t.b)=\tau_a(\tilde b(tua))=\tilde b(\tau(tua))
=\tilde b(t^{-1}\tau(ua))=t^{-1}.\tau_a(b), $$
establishing 2.
\par\nin 3. Write $N_+$ for the complex 
subgroup of $N_1$ corresponding to the Lie 
algebra 
$\n_+=\bigoplus_{\alpha\in \Sigma_1^+\atop 
\alpha\res_\a=0}\g_\C^\alpha$. Notice that 
$N_1=N_\C \rtimes N_+$. Fix $g\in G$. 
Then $g$ can be uniquely expressed as  
$$g=nn_+ a t u\leqno(4.3.1)$$
with $n\in N_\C$, $n_+\in N_+$, $a\in A$, 
$t\in \exp(i\t_1)$ and $u\in U$. 
Replacing $g$ with $\tau(g)$ one obtains 
a decomposition 
$$\tau(g)= n' n_+' a' t' u' \leqno(4.3.2)$$
for  $n'\in N_\C$, $n_+'\in N_+$, $a'\in A$, 
$t'\in \exp(i\t_1)$ and $u'\in U$. 
We claim that $a=a'$. Clearly this 
will prove the assertion in 3. 
\par We apply $\tau$ to (4.3.1)
$$\tau(g)=\tau(n)\tau(n_+) a t^{-1} \tau(u).$$
With (4.3.2) we now get the equality 
$$n' n_+' a' t' u' = \tau(n) \tau(n_+) a t^{-1} \tau(u).
\leqno(4.3.3)$$
Observe that $\tau(N_\C)=N_\C$ and $\tau(U)=U$, but 
$\tau(N_+)=N_-$. 
Let $\theta$ denote the Cartan-involution on $G_\C$ 
with fixed point group $U$. Symmetrizing (4.3.3) 
we obtain that 
$$\eqalign{\tau(g) \theta (\tau(g^{-1}))&=
n'n_+' (a')^2 (t')^2 \theta(n_+')^{-1} \theta(n')^{-1}\cr 
&=\tau(n)\tau(n_+) a^2 t^{-2} \theta\tau(n_+)^{-1} \theta\tau(n)^{-1}.}$$ 
Notice that $n_+' (a')^2 (t')^2 \theta(n_+')^{-1}$ and 
$\tau(n_+) a^2 t^{-2} \theta\tau(n_+)^{-1}$ belong to the reductive 
group $Z_{G_\C}(A)$. Thus the Bruhat-decomposition 
of $G_\C$ with respect to the parabolic subgroup 
$Z_{G_\C}(A) N_\C$ implies that $n'=\tau(n)$. 
Hence we obtain the identity  
$$n_+' (a')^2 (t')^2 \theta(n_+')^{-1}
=\tau(n_+) a^2 t^{-2} \theta\tau(n_+)^{-1} \leqno(4.3.4)$$ 
in $Z_{G_\C}(A)$. Notice that $A$ lies in the center of 
$Z_{G_\C}(A)$. Thus (4.3.4) implies that 
$$(a')^2 x'=a^2x\quad \hbox{ for some $x',x\in [Z_{G_{\C}}(A),Z_{G_{\C}}(A)]_0$}.  $$ 
Hence $(a')^2=a^2$ and so $a=a'$ as was to be shown. \qed  
\enddemo

\remark{Remark 4.3.2} (a) If $\m=\t_1$ is abelian, 
then the involution $\tau_a\: M_a\to M_a$ is 
anti-symplectic (cf.\ [8]). In fact, if $\m$ is abelian, then 
$\tau(B_1)=B_1$ and $\pr_\uu\circ \tau=\tau\circ \pr_\uu$ holds. 
One verifies the anti-symplecticity of $\tau_a$ along the 
same lines as Lemma 4.3.1(1). 
\par\nin (b) If $\m$ is not abelian, then 
$\tau(B_1)\neq B_1$ and $\pr_\uu\circ \tau\neq\tau\circ \pr_\uu$. 
One can verify from (4.2.2) that $\tau_a$ is not anti-symplectic in this 
case (see Example 4.3.3 below). Mistakenly, as pointed out in [8], it was 
implicitly assumed in [13] that $\tau_a$ is always anti-symplectic. \qed 
\endremark 

\msk\nin{\bf Proof of Theorem 4.1.1:}
In view of Lemma 4.3.1 the assumptions of 
Theorem 3.1 are satisfied and we can conclude 
that $\Phi(Q_a)=\Phi(M_a)=\conv (\Fix(M_a))$ . 
Standard structure theory shows that $\Fix(M_a)=
{\Cal W}.a$. 
\qed 

\example{Example 4.3.3} We will show that $\tau_a$ is 
not anti-symplectic in case $\m$ is not abelian. In this 
situation $\g$ contains 
a subalgebra of type $\so(1,4)$. Therefore it is 
enough to consider the case of $G=\SO_e(1,4)$. 

\par We start with some comments of general nature.  First notice 
that every $n\in N$ is contained in some $M_a=\tilde b(Ka)$ for 
some $a\in A$. Now fix $n\in N$ and $a\in A$ 
such that $n\in M_a$. Let $X,Y\in i\p$. Then 
$\tau_a(n)=n$ and 
$$d\tau_a(n)\tilde X_n =-\tilde X_n \quad\hbox{and}\quad 
d\tau_a(n)\tilde Y_n =-\tilde Y_n .$$
Thus if $\tau_a$ were anti-symplectic, then $\omega_n(\tilde X_n, \tilde Y_n)=0$. 
Below we will show that $\omega_n(\tilde X_n, \tilde Y_n)\neq 0$ for 
a specific choice of elements $n, X,Y$. 
\ssk Let now $G=\SO_e(1,4)$. Then the Lie algebra of $G$ is 
given by
$$\g=\left\{\pmatrix 0 & u^t \\ u & X\endpmatrix\: u\in \R^4, \ X\in \so(4)\right\}.$$
The complexification of $\g$ then is 
$$\g_\C=\left\{\pmatrix 0 & w^t \\ w & Z\endpmatrix\: w\in \C^4, \ Z\in \so(4,\C)\right\}.$$
Our choice of $\a$ and $\t_1$ will be 

$$\a=\R\pmatrix  0 & 0 & 0& 0 & 1 \\ 0 & 0 & 0& 0 & 0 \\ 0 & 0 & 0& 0 & 0 \\
 0 & 0 & 0& 0 & 0 \\ 1 & 0 & 0& 0 & 0 \endpmatrix 
\quad\hbox{and}\quad \t_1=\R\pmatrix  0 & 0 & 0& 0 & 0 \\ 0 & 0 & 1& 0 & 0 \\ 0 & -1 & 0& 0 & 0 \\
 0 & 0 & 0& 0 & 0 \\ 0 & 0 & 0& 0 & 0 \endpmatrix .$$
With appropriate $\Sigma_1^+$ the nilpotent Lie algebras $\n_\C$ and $\n_+$
are given by 
$$\n_\C=\left\{\pmatrix  0 & z^t & 0\\ z& 0&-z\\ 0& z^t &0 \endpmatrix \: z\in \C^3\right\}\ , $$
and 
$$\n_1=\C \pmatrix  0 & 0 & 0& 0 & 0 \\ 0 & 0 & 0& i & 0 \\ 0 & 0 & 0& 1 & 0 \\
 0 & -i & -1& 0 & 0 \\ 0 & 0 & 0& 0 & 0 \endpmatrix.$$

\ssk Define an element $n\in N$ by 
$$n=\pmatrix  {5\over 2} & 1 & 1& 1& -{3\over 2} \\ 1 & 1 & 0& 0 & -1 \\ 1 & 0 & 1& 0 & -1 \\
1 & 0 & 0& 1 & -1 \\ {3\over 2} & 1 & 1& 1 & -{1\over 2} \endpmatrix\in N.$$
Next we introduce elements $X,Y\in i\p$ by 

$$X=\pmatrix  0 & 0 & 0& i & 0 \\ 0 & 0 & 0& 0 & 0 \\ 0 & 0 & 0& 0 & 0 \\
i & 0 & 0& 0 & 0 \\ 0 & 0 & 0& 0 & 0 \endpmatrix
\quad \hbox{and}\quad Y=\pmatrix  0 & i & 0& 0 & 0 \\ i & 0 & 0& 0 & 0 \\ 0 & 0 & 0& 0 & 0 \\
0 & 0 & 0& 0 & 0 \\ 0 & 0 & 0& 0 & 0 \endpmatrix$$
A simple computation gives 

$$\Ad(n)^{-1} X=\pmatrix  0 & -i & -i& {3\over 2}i  & -i \\ -i & 0 & 0& -i & i \\ -i& 0 & 0& -i & i \\
{3\over 2}i & i & i& 0 & -{1\over 2}i\\ -i& -i & -i& {1\over 2}i & 0 \endpmatrix$$
and 
$$\Ad(n)^{-1} Y=\pmatrix  0 & {3\over 2}i & -i& -i & -i \\ {3\over 2}i & 0 & i& i & -{1\over 2}i \\ -i& -i & 0& 0 & i \\
-i & -i & 0& 0 & i\\ -i& {1\over 2}i & -i& -i & 0 \endpmatrix.$$
Then  
$$\pr_\uu (\Ad(n)^{-1}X)=\pmatrix  0 & 0 & 0& i & -i \\ 0 & 0 & 0& -1 & 0 
\\ 0 & 0 & 0& 1 & 0 \\
i & 1 & -1& 0 & 0 \\ -i & 0 & 0& 0 & 0 \endpmatrix.$$
Up to scalar we may identify $\B$ with the imaginary 
part of the trace form. With this normalization we then have 

$$\eqalign{\omega_n(\tilde X_n, \tilde Y_n)&=\Im \tr\left(
\pmatrix  0 & 0 & 0& i & -i \\ 0 & 0 & 0& -1 & 0 
\\ 0 & 0 & 0& 1 & 0 \\
i & 1 & -1& 0 & 0 \\ -i & 0 & 0& 0 & 0 \endpmatrix\
\cdot \pmatrix  0 & {3\over 2}i & -i& -i & -i \\ {3\over 2}i & 0 & i& i & -{1\over 2}i \\ -i& -i & 0& 0 & i \\
-i & -i & 0& 0 & i\\ -i& {1\over 2}i & -i& -i & 0 \endpmatrix\right)\cr
&=\Im \tr\pmatrix  0 & * & *& * & * \\ * & i & *& * & * 
\\ * & * & 0& *& * \\
* & * & *& 1+i & * \\ * & * & *& * & -1 \endpmatrix=2\neq 0.\cr}$$
\endexample

Following a suggestion from the referee we want to add a few remarks.
\remark{Remark 4.3.4} (a) Theorem 4.1.1 is the "real nonlinear" version of Kostant's convexity theorem. To complete its symplectic proof suggested in [13] we needed our generalized symplectic convexity theorem 3.1. There is also a (real) linear version of Kostant's theorem:
Let $Y\in \a$. Then $$ a(\Ad(K).Y)=\conv ({\Cal W}.Y), $$
where $\conv(\cdot)$ denotes the convex hull of $(\cdot)$ and $a:\g=\k+\a+\n \rightarrow \a$ denotes the middle projection for the Iwasawa decomposition on the Lie algebra level. \par\nin
In [4] and [13] symplectic proofs for this result have been given both for complex and for real $G$. In this case Duistermaat's theorem (Theorem 1.2 in the introduction) suffices, because complex conjugation restricted to the coadjoint orbit $\Ad(K).Y$ is antisymplectic with respect to the natural symplectic form on $\Ad(K).Y$ for any complex semisimple group $G$.
\par (b) The linearization trick of Alekseev-Meinrenken-Woodward [1] (Chapter 3) can be generalized as follows. We adopt notation and numbering from [1]. \par\nin
Let $K$ be a compact connected Lie group equipped with an involutive automorphism $\sigma_K$ of the type described in Section 2.4 in [1]. Suppose $M$ is a $K$-manifold with an involution $\sigma_M$. Let $\Omega, \omega$ be two forms on $M$, and $\Psi:M\rightarrow K^*, \ \Phi:M\rightarrow \k^*$ related by (10) (i.e. equation (10) in [1]).
Denote by $Q$ the fixed point set of $M$ under $\sigma_M$. Then the following are equivalent:
\roster
\item $(M,\Omega,\Psi)$ is a Hamiltonian $K$-space with $K^*$-valued moment map, and the involution $\sigma_M$ satisfies \par
(i) $\sigma_M(k.m)=\sigma_K(k).(\sigma_M(m)) \quad \forall \ k \in K, m\in M $ \par
(ii) $\Psi \circ \sigma_M = \sigma_{K^*} \circ \Psi$ \par
(iii) $Q$ is a Lagrangian subspace of the symplectic manifold $(M,\Omega)$.
\item $(M,\omega,\Phi)$ is a Hamiltonian $K$-space with $\k^*$-valued moment map, and the involution $\sigma_M$ satisfies \par
(i) $\sigma_M(k.m)=\sigma_K(k).(\sigma_M(m)) \quad \forall \ k\in K, m\in M $ \par
(ii) $\Phi \circ \sigma_M = \sigma_{\k^*} \circ \Phi$ \par
(iii) $Q$ is a Lagrangian subspace of the symplectic manifold $(M,\omega)$.
\endroster
This follows easily from (15) and (16). \par
Note that (i), (ii), (iii) become our properties (2.1.2), (2.1.3), (2.1.4) in the case $K=T$ and $\sigma_K(k)=k^{-1}$. \par
Whereas linearization as in 3.1 in [1] can be generalized to a set-up with non anti-symplectic involutions, this is not clear for the results of 3.3 (and the corresponding statements in 3.4), since their proofs depend on the invariance of the flow $\tilde{\phi}^s$ under $\sigma_M$ in the isotopy lemma 3.4. 
\par (c) Duistermaat's generalized result on the tightness of Hamiltonian functions and its proof (Section 3 in [4]) also hold under the weaker conditions (2.1.2), (2.1.3), (2.1.4) on the involution $\tau$.

\endremark

\refstyle{C}

\enddocument